\documentclass[12pt]{article}
\usepackage{amsmath}

\usepackage{t1enc}
\usepackage[latin1]{inputenc}
\usepackage[english]{babel}
\usepackage{amsmath,amsthm}
\usepackage{amsfonts}
\usepackage{latexsym}
\usepackage[dvips]{graphicx}
\usepackage{graphicx}
\usepackage[natural]{xcolor}
\newtheorem{theorem}{Theorem}
\newtheorem{lemma}[theorem]{Lemma}
\newtheorem{proposition}[theorem]{Proposition}

\newtheorem{example}[theorem]{Example}
\newtheorem{definition}[theorem]{Definition}

\begin{document}

\date{}
\author{Liangquan Zhang \ and \ \ Yufeng SHI \thanks{
Supported by National Natural Science Foundation of China Grant 10771122,
Natural Science Foundation of Shandong Province of China Grant Y2006A08 and
National Basic Research Program of China (973 Program, No. 2007CB814906).}
\thanks{
Corresponding author, E-mail: yfshi@sdu.edu.cn, Tel: +86-531-88364100, Fax:
+86-531-88564100.} \\
School of Mathematics, Shandong University\\
Jinan 250100, People's Republic of China}
\title{\textbf{Comparison Theorem of Multi-dimensional Backward Doubly Stochastic
Differential Equations on Infinite Horizon}}
\maketitle

\begin{abstract}
Under quasi-monotone assumptions for coefficients, we show one kind of
comparison theorem for multi-dimensional\textbf{\ }backward doubly
stochastic differential equations on infinite horizon. An example is given
as well.
\end{abstract}



\section{Introduction}

\hspace{0.2 in}Since the fundamental paper on nonlinear backward stochastic
differential equations (BSDEs in short) was published by Pardoux and Peng$%
^{[1]}$ in 1990, BSDEs have attracted greatly interests from both
mathematics community (e.g. see [2] and its references) and finance
community (e.g. see [3] and its references). In the theory of BSDEs, the
comparison theorems are very important results and powerful tools to treat
many problems such as stochastic differential games, $g$-expectations,
viscosity solutions of partial differential equations, options prices and so
on (for instance, see [3] for the comparison theorem of one-dimensional
BSDEs and its applications). In 1994, Christel and Ralf$^{\left[ 4\right] }$
proved the comparison theorems for finite and infinite dimensional
stochastic differential equations (SDEs in short). Enlightened from this
idea, Zhou$^{\left[ 5\right] }$ obtained a comparison theorem for the
multi-dimensional BSDEs in finite time intervals in 1999.

A class of backward doubly stochastic differential equations (BDSDEs in
short) was introduced by Pardoux and Peng$^{\left[ 6\right] }$ in 1994, in
order to provide a probabilistic interpretation for the solutions to a class
of quasilinear stochastic partial differential equations (SPDEs in short).
BDSDEs are a new kind of BSDEs with two different directions of stochastic
integrals, i.e. the equations involve both a standard (forward) stochastic
integral d$\overrightarrow{W_t}$ and a backward stochastic integral d$%
\overleftarrow{B_t}$. In [6] the authors have proved the existence and
uniqueness of solutions to BDSDEs under uniformly Lipschitz assumptions on
coefficients $f,$ $g$, i.e. for any square integrable terminal value $\xi ,$
the following BDSDE has a unique solution pair $\left( y_t,z_t\right) $ in
the interval $\left[ 0,T\right] .$%
\begin{equation}
y_t=\xi +\int_t^Tf\left( s,y_s,z_s\right) \text{d}s+\int_t^Tg\left(
s,y_s,z_s\right) \text{d}\overleftarrow{B_s}-\int_t^Tz_s\text{d}%
\overrightarrow{W_s},\text{ }t\in \left[ 0.T\right] .  \tag{1.1}
\end{equation}
Pardoux and Peng$^{\left[ 6\right] }$ have showed the relationship between
BDSDEs and the quasilinear SPDE as follows:
\[
\left\{
\begin{array}{c}
u\left( t,x\right) =h\left( x\right) +\int_t^T\left[ \mathcal{L}u\left(
s,x\right) +f\left( s,x,u\left( s,x\right) ,\left( \nabla u\sigma \right)
\left( s,x\right) \right) \right] \text{d}s \\
+\int_t^Tg\left( s,x,u\left( s,x\right) ,\left( \nabla u\sigma \right)
\left( s,x\right) \right) \text{d}\overleftarrow{B_s},\quad 0\leq t\leq T,
\end{array}
\right.
\]
where $u:\left[ 0,T\right] \times \mathbf{R}^k\rightarrow \mathbf{R}^k$ and $%
\nabla u\left( s,x\right) $ denotes the first order derivative of $u\left(
s,x\right) $ with respect to $x$, and
\[
\mathcal{L}u=\left(
\begin{array}{c}
Lu_1 \\
\vdots \\
Lu_k
\end{array}
\right) ,
\]
with $L\phi \left( x\right) =\frac 12\sum_{i,j=1}^d\left( \sigma \sigma
^{*}\right) _{ij}\left( x\right) \frac{\partial ^2\phi \left( x\right) }{%
\partial x_i\partial x_j}+\sum_{i=1}^db_i\left( x\right) \frac{\partial \phi
\left( x\right) }{\partial x_i}$ (for more details see [6]).

As a result of their important significance to the theory of stochastic
partial differential equations, BDSDEs have attracted more and more
researchers' interest, for example, Bally and Matoussi$^{\left[ 7\right] }$,
Shi, Gu and Liu$^{\left[ 8\right] }$, Han, Zhu and Shi$^{\left[ 9\right] }$,
Zhang and Zhao$^{\left[ 10\right] }$, Ren, Lin and Hu$^{\left[ 11\right] }$
and their references. Furthermore, by means of the probabilistic
representations of SPDEs, it is available to provide numerical simulations
of SPDEs by BDSDEs (cf. Shi, Yang and Yuan$^{\left[ 12\right] }$). So far,
however most of works involve BDSDEs in finite time intervals. There is a
little work to treat BDSDEs for infinite time interval cases. Infinite
horizon BDSDEs are very important to provide a probabilistic representation
for some quasilinear stochastic partial differential equations. For
instance, Zhang and Zhao$^{\left[ 10\right] }$ proved the existence and
uniqueness of solutions for BDSDEs on infinite horizons, and described the
stationary solutions of some SPDEs by virtue of the solutions of BDSDEs on
infinite horizons. Recently, Han, Shi and Zhu$^{\left[ 13\right] }$ applied
the arguments introduced by Chen and Wang$^{\left[ 14\right] }$ to get the
existence and uniqueness of solutions of the following BDSDE:
\begin{equation}
y_t=\xi +\int_t^\infty f\left( s,y_s,z_s\right) \text{d}s+\int_t^\infty
g\left( s,y_s\right) \text{d}\overleftarrow{B_s}-\int_t^\infty z_s\text{d}%
\overrightarrow{W_s},\text{ }t\in \left[ 0,\infty \right) .  \tag{1.2}
\end{equation}

The aim of this paper is to study further the theory of infinite horizon
BDSDE (1.2). In this paper, under quasi-monotonicity conditions for
coefficient $f$, we show a comparison theorem for multi-dimensional BDSDE
(1.2).

This paper is organized as follows: We present the setting of problems, the
preliminaries and assumptions in Section 2. In Section 3 we prove a
comparison theorem for multi-dimensional BDSDE (1.2) on infinite horizon, by
an argument similar to Zhou$^{\left[ 5\right] }$. Lastly, we give an example
to show our comparison theorem of BDSDEs.

\section{Preliminaries: the existence and uniqueness of infinite horizon
BDSDEs and an extension of It\^o's formula}

\hspace{0.2 in}\textbf{Notation}. We denote the Euclidean norm of a vector $%
x\in R^k$ by $|x|,$ and the one of a $d\times k$ matrix $A$ by $\parallel
A\parallel =\sqrt{TrAA^{*}},$ where $A^{*}$ is the transpose of $A$.

Let $(\Omega ,\mathcal{F},P)$ be a complete probability space, $%
\{W_t\}_{t\geq 0}$ and $\{B_t\}_{t\geq 0}$ be two mutually independent
standard Brownian motions, with values respectively in $R^d$ and in $R^l$,
defined on $(\Omega ,\mathcal{F},P)$. Let $\mathcal{N}$ denote the class of $%
P$-null sets of $\mathcal{F}$. For each $t\in \left[ 0,\infty \right] ,$ we
define
\begin{eqnarray*}
\mathcal{F}_{0,t}^W &\doteq &\sigma \left\{ W_r;0\leq r\leq
t\right\} \vee \mathcal{N},\quad \ \mathcal{F}_{t,\infty }^B\doteq
\sigma \left\{
B_r-B_t;t\leq r<\infty \right\} \vee \mathcal{N}\text{,} \\
\mathcal{F}_{0,\infty }^W &\doteq &\bigvee\limits_{0\leq t<\infty }\mathcal{F%
}_{0,t}^W,\quad \ \mathcal{F}_{\infty ,\infty }^B\doteq
\bigcap\limits_{0\leq t<\infty }\mathcal{F}_{t,\infty }^B\text{,}
\end{eqnarray*}

and
\[
\mathcal{F}_t\doteq \mathcal{F}_{0,t}^W\vee \mathcal{F}_{t,\infty }^B,\quad
t\in \left[ 0,\infty \right] .
\]
Note that $\left\{ \mathcal{F}_{0,t}^W;t\in \left[ 0,\infty \right] \right\}
$ is increasing and $\left\{ \mathcal{F}_{t,\infty }^B;t\in \left[ 0,\infty
\right] \right\} $ is decreasing, and the collection $\left\{ \mathcal{F}%
_t;t\in \left[ 0,\infty \right] \right\} $ is neither increasing nor
decreasing. Suppose
\[
\mathcal{F}=\mathcal{F}_{0,\infty }^W\vee \mathcal{F}_{0,\infty }^B.
\]

For any $n\in N,$ let $\mathcal{S}^2(\mathbf{R}^{+};\mathbf{R}^n)$ denote
the space of all $\{\mathcal{F}_t\}$-measurable $n$-dimensional processes $%
\upsilon $ with norm of $\parallel \upsilon \parallel _{\mathcal{S}}\doteq \left[ \mathbf{E}%
(\sup\limits_{s\geq 0}|\upsilon (s)|)^2\right] ^{\frac 12}<\infty .$

\begin{definition}
A stochastic process $X=\left\{ X_t;t\geq 0\right\} $ is called $\mathcal{F}%
_t$-progressively measurable, if for any $0\leq t\leq T\leq \infty $, $X$ on
$\Omega \times \left[ 0,t\right] $ is measurable with respect to $\left(
\mathcal{F}_{0,t}^W\times \mathcal{B}\left( \left[ 0,t\right] \right)
\right) \vee \left( \mathcal{F}_{t,\infty }^B\times \mathcal{B}\left( \left[
t,T\right] \right) \right) $.
\end{definition}

We denote similarly by $\mathcal{M}^2\left( \mathbf{R}^{+};\mathbf{R}%
^n\right) $ the space of (class of d$P\otimes $d$t$ a.e equal) all $\left\{
\mathcal{F}_t\right\} $-measurable $n$-dimensional processes $\upsilon $
with norm of $\parallel \upsilon \parallel _{\mathcal{M}}\doteq \left[
\mathbf{E}\int_0^\infty |\upsilon (s)|^2\text{d}s\right] ^{\frac 12}<\infty .
$

For any $t\in \left[ 0,\infty \right] ,$ let $L^2\left( \Omega ,\mathcal{F}%
_t,P;\mathbf{R}^n\right) $ denote the space of all $\left\{ \mathcal{F}%
_t\right\} $-measurable $\mathbf{R}^n$-valued random variable $\xi $
satisfying $\mathbf{E}\left| \xi \right| ^2<\infty .$ We also denote that
\[
\mathcal{B}^2\doteq \left\{ \left( X,Y\right) :X\in \mathcal{S}^2\left(
\mathbf{R}^{+};\mathbf{R}^n\right) ,\quad Y\in \mathcal{M}^2\left( \mathbf{R}%
^{+};\mathbf{R}^n\right) \right\} .
\]
For each $(X,Y)\in \mathcal{B}^2$, we define the norm of $(X,Y)$ by
\[
\left\| \left( X,Y\right) \right\| _{\mathcal{B}}=\left( \left\| X\right\| _{%
\mathcal{S}}^2+\left\| Y\right\| _{\mathcal{M}}^2\right) ^{\frac 12}.
\]
Obviously $\mathcal{B}^2$ is a Banach space.

Consider the following infinite horizon backward doubly stochastic
differential equations
\begin{equation}
y_t=\xi +\int_t^\infty f\left( s,y_s,z_s\right) \text{d}s+\int_t^\infty
g\left( s,y_s\right) \text{d}\overleftarrow{B_s}-\int_t^\infty z_s\text{d}%
\overrightarrow{W_s},\text{ }0\leq t\leq \infty ,  \tag{2.1}
\end{equation}
where $\xi \in L^2\left( \Omega ,\mathcal{F}_\infty ,P;\mathbf{R}^k\right) $
is given. We note that the integral with respect to $\left\{ B_t\right\} $
is a ``backward It\^o integral'' and the integral with respect to $\left\{
W_t\right\} $ is a standard forward It\^o integral. These two types of
integrals are particular cases of the It\^o-Skorohod integral, see Nualart
and Pardoux$^{\left[ 15\right] }$.

All the equalities and inequalities mentioned in this paper are in sense of d%
$t\times $d$p$ almost surely on $\left[ 0,\infty \right] \times
\Omega .$ In addition, throughout this paper, when we indicate two
vectors $\nu ^1,$ $\nu ^2$ satisfying $\nu ^1\geq \nu ^2$, it means
that we have $\nu ^{1,j}\geq \nu ^{2,j}$ for each $j$ of their
components.

\begin{definition}
A pair of processes $\left( y,z\right) $: $\Omega \times \mathbf{R}%
^{+}\rightarrow \mathbf{R}^k\times \mathbf{R}^{k\times d}$ is called a
solution of BDSDE (2.1), if $\left( y,z\right) \in \mathcal{B}^2$ and
satisfies BDSDE (2.1) for any $t\in \left[ 0,\infty \right] $.
\end{definition}

Let
\begin{eqnarray*}
f &:&\Omega \times \mathbf{R}^{+}\times \mathbf{R}^k\times \mathbf{R}%
^{k\times d}\rightarrow \mathbf{R}^k, \\
g &:&\Omega \times \mathbf{R}^{+}\times \mathbf{R}^k\rightarrow \mathbf{R}%
^{k\times l},
\end{eqnarray*}
satisfy the following assumptions:

\begin{enumerate}
\item[(H1)]  for any $(y,z)\in \mathbf{R}^k\times \mathbf{R}^{k\times d},$ $%
f\left( \cdot ,y,z\right) $ and $g\left( \cdot ,y\right) $ are $\left\{
\mathcal{F}_t\right\} $-progressively measurable process, such that
\[
\mathbf{E}\left( \int_0^\infty f(t,0,0)\text{d}t\right) ^2<\infty ;\text{ }%
g\left( \cdot ,0\right) \in \mathcal{M}^2\left( \mathbf{R}^{+};\mathbf{R}%
^{k\times l}\right) .
\]

\item[(H2)]  $f$ and $g$ satisfy Lipschitz conditions with Lipschitz
functions $\upsilon :=\left\{ \upsilon (t)\right\} $ and $u:=\left\{
u(t)\right\} ,$ that is, there exists two positive deterministic function $%
\left\{ \upsilon \left( t\right) \right\} $ and $\left\{ u\left( t\right)
\right\} $ such that:
\[
\left\{
\begin{array}{l}
|f\left( t,y_1,z_1\right) -f\left( t,y_2,z_2\right) |\leq \upsilon \left(
t\right) \left| y_1-y_2\right| +u\left( t\right) \left| z_1-z_2\right| , \\
\left| g\left( t,y_1,z_1\right) -g\left( t,y_2,z_2\right) \right| \leq
u\left( t\right) \left| y_1-y_2\right| , \\
\forall \left( t,y_i,z_i\right) \in \mathbf{R}^{+}\times \mathbf{R}^k\times
\mathbf{R}^{k\times d},\quad i=1,2.
\end{array}
\right.
\]
\end{enumerate}

For the existence and uniqueness of solutions to BDSDE (2.1), we need the
following conditions on $v$ and $u$:

\begin{enumerate}
\item[(H3)]  $\int_0^\infty v(t)$d$t<\infty $ and $\int_0^\infty u^2(t)$d$t<\infty $.
Furthermore $v\left( t\right) $ and $u\left( t\right) $ either are
bounded on R$^{+}$, or are unbounded but at most have finite many
discontinuous points.
\end{enumerate}

\begin{proposition}
Let $\xi \in L^2\left( \Omega ,\mathcal{F}_\infty ,P;\mathbf{R}^k\right) $
be given, (H1), (H2) and (H3) hold for $f$ and $g$. Then the BDSDE (2.1) has
a unique solution $\left( y,z\right) \in \mathcal{B}^2.$
\end{proposition}

This proposition was proved in [13]. We shall employ the following extension
of It\^o's formula.

\begin{proposition}
Let $\alpha \in \mathcal{S}^2\left( [0,\infty ];\mathbf{R}^k\right) ,$ $%
\beta \in \mathcal{M}^2\left( [0,\infty ];\mathbf{R}^k\right) ,$ $\gamma \in
\mathcal{M}^2\left( [0,\infty ];\mathbf{R}^{k\times l}\right) ,$ $\delta \in
\mathcal{M}^2\left( [0,\infty ];\mathbf{R}^{k\times d}\right) $ satisfy:
\[
\alpha _t=\alpha _0+\int_0^t\beta _s\text{d}s+\int_0^t\gamma _s\text{d}%
\overleftarrow{B_s}+\int_0^t\delta _s\text{d}\overrightarrow{W_s}\text{%
,\qquad }0\leq t\leq T.
\]

\noindent Then
\begin{eqnarray*}
\left| \alpha _t\right| ^2=\left| \alpha _0\right| ^2+2\int_0^t\left( \alpha
_s,\beta _s\right) \text{d}s+2\int_0^t\left( \alpha _s,\gamma _s\text{d}%
\overleftarrow{B_s}\right) +2\int_0^t\left( \alpha _s,\delta _s\text{d}%
\overrightarrow{W_s}\right)  \\
-\int_0^t\parallel \gamma _s\parallel ^2\text{d}s+\int_0^t\parallel \delta
_s\parallel ^2\text{d}s, \\
\mathbf{E}\left| \alpha _t\right| ^2=\mathbf{E}\left| \alpha _0\right| ^2+2%
\mathbf{E}\int_0^t\left( \alpha _s,\beta _s\right) \text{d}s-\mathbf{E}%
\int_0^t\parallel \gamma _s\parallel ^2\text{d}s+\mathbf{E}\int_0^t\parallel
\delta _s\parallel ^2\text{d}s.
\end{eqnarray*}
\end{proposition}

This proposition was also obtained in [6].

\section{Comparison theorem of multi-dimensional BDSDEs on infinite horizon}

\hspace{0.2 in}Now we consider the following $k$-dimensional BDSDEs on
infinite horizon:
\begin{equation}
y_t^1=\xi ^1+\int_t^\infty f^1\left( s,y_s,z_s\right) \text{d}%
s+\int_t^\infty g\left( s,y_s\right) \text{d}\overleftarrow{B_s}%
-\int_t^\infty z_s^1\text{d}\overrightarrow{W_s},\text{ }0\leq t\leq \infty ,
\tag{3.1}
\end{equation}
\begin{equation}
y_t^2=\xi ^2+\int_t^\infty f^2\left( s,y_s,z_s\right) \text{d}%
s+\int_t^\infty g\left( s,y_s\right) \text{d}\overleftarrow{B_s}%
-\int_t^\infty z_s^2\text{d}\overrightarrow{W_s},\text{ }0\leq t\leq \infty ,
\tag{3.2}
\end{equation}
where $\xi ^1,$ $\xi ^2\in L^2\left( \Omega ,\mathcal{F}_\infty ,P;\mathbf{R}%
^k\right) $%
\begin{eqnarray*}
f^1\left( \omega ,t,y,z\right) ,f\left( \omega ,t,y,z\right) &:&\Omega
\times \left[ 0,\infty \right] \times \mathbf{R}^k\times \mathbf{R}^{k\times
d}\rightarrow \mathbf{R}^k, \\
g\left( \omega ,t,y\right) &:&\Omega \times \left[ 0,\infty \right] \times
\mathbf{R}^k\rightarrow \mathbf{R}^{k\times l}.
\end{eqnarray*}

We assume that $\xi ^1,$ $\xi ^2$ and $f^1,$ $f^2,$ $g$ satisfy the
assumptions (H1)-(H3), in addition to the following conditions:

\begin{enumerate}
\item[(H4)]
\[
\left\{
\begin{array}{l}
\text{(i)\quad }\xi ^1\geq \xi ^2; \\
\text{(ii)\quad For all }j=1,2,\cdots k,\text{ }\forall \omega \times t\in
\Omega \times \left[ 0,\infty \right] ,\text{ }f_j^1\left( \omega
,t,y^1,z^1\right) \geq f_j^2\left( \omega ,t,y^2,z^2\right) , \\
\text{\qquad where }y^1,\text{ }y^2\in \mathbf{R}^k,\text{ }z^1,\text{ }%
z^2\in \mathbf{R}^{k\times d},\text{ }y_j^1=y_j^2,\text{ }z_j^1=z_j^2,\text{
}y_l^1\geq y_l^2,\text{ }l\neq j;
\end{array}
\right.
\]
\end{enumerate}

From Proposition 3, under the assumptions (H1)-(H3), there exist two pairs
of measurable processes $\left( y^1,z^1\right) \in \mathcal{S}^2\left(
\left[ 0,\infty \right] ;\mathbf{R}^k\right) \times \mathcal{M}^2\left(
0,\infty ;\mathbf{R}^{k\times d}\right) $ and $\left( y^2,z^2\right) \in
\mathcal{S}^2\left( [0,\infty ];\mathbf{R}^k\right) \times \mathcal{M}%
^2\left( 0,\infty ;\mathbf{R}^{k\times d}\right) $ satisfying BDSDE (3.1)
and (3.2), respectively. Before showing our major result, we need to prove a
kind of Gronwall's lemma.

\begin{lemma}
Suppose $m\left( t\right) ,$ $r\left( t\right) \in C\left( \mathbf{R}^{+},%
\mathbf{R}^{+}\right) ,$ $\int_0^\infty r\left( s\right) $d$s<\infty .$ If
\begin{equation}
m\left( t\right) \leq A+M\int_0^\infty m\left( s\right) r\left( s\right)
\text{d}s,\text{ }t\in \mathbf{R}^{+},  \tag{3.3}
\end{equation}
\noindent where $A\geq 0,$ $M>0,$ then we have
\[
m\left( t\right) \leq A\exp \left[ M\int_t^\infty r\left( s\right) \text{d}%
s\right] ,\text{ }\forall t\in \mathbf{R}^{+}.
\]
\end{lemma}

\noindent \textbf{Proof.}\quad For any $T\in \mathbf{R}^{+},$ let us
set
\[
H\left( t\right) =\int_t^Tm\left( s\right) r\left( s\right) \text{d}s,\text{
}g\left( t\right) =H\left( t\right) \exp \left\{ -M\int_t^Tr\left( s\right)
\text{d}s\right\} ,\quad t\in \left[ 0,T\right] .
\]
It is easy to check that $g\left( T\right) =0,$ and $H^{^{\prime }}\left(
t\right) =-m\left( t\right) r\left( t\right) .$ Then
\begin{equation}
g^{^{\prime }}\left( t\right) =H^{^{\prime }}\left( t\right) \exp \left\{
M-\int_t^Tr\left( s\right) \text{d}s\right\} +MH\left( t\right) r\left(
t\right) \exp \left\{ -M\int_t^Tr\left( s\right) \text{d}s\right\} .
\tag{3.4}
\end{equation}
From (3.4), we have
\begin{equation}
g^{^{\prime }}\left( t\right) \geq -Ar\left( t\right) \exp \left\{
M-\int_t^Tr\left( s\right) \text{d}s\right\} .  \tag{3.5}
\end{equation}
Integrating [3.5] on $\left[ t,T\right] ,$ we obtain
\[
g\left( t\right) \leq A\int_t^Tr\left( s\right) \exp \left\{
-M\int_s^Tr\left( \nu \right) d\nu \right\} \text{d}s.
\]
That is,
\[
H\left( t\right) \leq A\int_t^Tr\left( s\right) \exp \left\{
M\int_t^sr\left( \nu \right) d\nu \right\} \text{d}s.
\]
Using (3.3) again, we have
\[
m\left( t\right) \leq A+MH\left( t\right) \leq A\exp \left\{
M\int_t^Tr\left( s\right) \text{d}s\right\} .
\]
Letting $T\rightarrow \infty ,$ we get
\[
m\left( t\right) \leq A\exp \left\{ M\int_t^\infty r\left( s\right) \text{d}%
s\right\} ,\text{ \quad }t\in R^{+}.
\]
The proof is completed.\quad $\Box$

\begin{theorem}
Assume the conditions (H1)-(H4) hold. Let $\left( y^1,z^1\right) $ and $%
\left( y^2,z^2\right) $ be the solutions of BDSDE (3.1) and (3.2),
respectively. Then $y^1\geq y^2,$ $\forall t\in \mathbf{R}^{+}.$
\end{theorem}

\noindent \textbf{Proof.}\quad $\forall \varepsilon >0,$ we define
the following
function $\phi \left( s\right) :\mathbf{R}\rightarrow \mathbf{R}$%
\[
\phi _\varepsilon \left( y\right) =\left\{
\begin{array}{l}
y^2,\text{\quad }y\leq 0, \\
y^2-\frac 1{6\varepsilon }y^3,\text{\quad }0\leq y\leq 2\varepsilon , \\
2\varepsilon y-\frac 43\varepsilon ^2,\text{\quad }y\geq 2\varepsilon .
\end{array}
\right.
\]
Clearly, $\phi _\varepsilon \left( y\right) \in C^2\left( \mathbf{R}\right) ,
$ and $\phi _\varepsilon ^{^{\prime \prime }}\left( y\right) $ is bound, for
all $y\in \mathbf{R},$%
\[
\phi _\varepsilon \left( y\right) \rightarrow \left| y^{-}\right| ^2,\quad
\text{ }\phi _\varepsilon ^{^{\prime }}\left( y\right) \rightarrow
-2y^{-},\quad \text{ }\phi _\varepsilon ^{^{\prime \prime }}\left( y\right)
\rightarrow 2\mathbf{I}_{\left\{ y<0\right\} },\text{ \quad as }\varepsilon
\rightarrow 0,
\]
where the first and second convergence are uniform, the last one is
convergence by point. We set
\[
\bar y_t=y_t^1-y_t^2,\quad \text{ }\bar z=z_t^1-z_t^2,\quad \text{ }\bar \xi
=\xi ^1-\xi ^2.
\]
Then $\left( \bar y_t,\bar z_t\right) $ satisfy the following BDSDE:
\[
\bar y=\bar \xi +\int_t^\infty \left[ f^1\left( s,y_s^1,z_s^1\right)
-f^2\left( t,y_s^2,z_s^2\right) \right] \text{d}s+\int_t^\infty \left[
g\left( s,y_s^1\right) -g\left( s,y_s^2\right) \right] \text{d}%
\overleftarrow{B_s}-\int_t^\infty \bar z_s\text{d}\overrightarrow{W_s}.
\]
For $\forall 0\leq t\leq T<\infty ,$ applying the extention of It\^o's
formula to $\phi _\varepsilon (\bar y_j(t))$ on $\left[ t,T\right] $, it
follows that
\begin{eqnarray*}
\phi _\varepsilon \left( \bar y_j\left( t\right) \right)  &=&\phi
_\varepsilon \left( \bar \xi _j\right) +\int_t^T\phi _\varepsilon ^{^{\prime
}}\left( \bar y_j\left( s\right) \right) \left[ f_j^1\left(
s,y_s^1,z_s^1\right) -f_j^2\left( s,y_s^2,z_s^2\right) \right] \text{d}s \\
&&\ \ +\int_t^T\phi _\varepsilon ^{^{\prime }}\left( \bar y_j(s)\right)
\left[ g_j\left( s,y_s^1\right) -g_j\left( s,y_s^2\right) \right] \text{d}%
\overleftarrow{B_s}-\int_t^T\phi _\varepsilon ^{^{\prime }}\left( \bar
y_j(s)\right) \bar z_s\text{d}\overrightarrow{W_s} \\
&&\ \ +\frac 12\int_t^T\phi _\varepsilon ^{^{^{\prime \prime }}}\left( \bar
y_j\left( s\right) \right) \left[ g_j\left( s,y_s^1\right) -g_j\left(
s,y_s^2\right) \right] ^2\text{d}s-\frac 12\int_t^T\phi _\varepsilon
^{^{^{\prime \prime }}}\left( \bar y_j\left( s\right) \right) \left| \bar
z_j(s)\right| ^2\text{d}s.
\end{eqnarray*}
Letting $T\rightarrow \infty $ and $\varepsilon \rightarrow 0,$
respectively, we get
\begin{eqnarray*}
\left| \left( \bar y_j\left( t\right) \right) ^{-}\right| ^2
&=&\left| \left( \bar \xi _j\right) ^{-}\right| ^2-\int_t^\infty
2\left( \bar y_j\left( s\right) \right) ^{-}\left[
f_j^1(s,y_s^1,z_s^1)-f_j^2(s,y_s^2,z_s^2)\right] \text{d}s \\
&&\ -2\int_t^\infty \left( \bar y_j\left( s\right) \right)
^{-}\left[ g_j(s,y_s^1)-g_j(s,y_s^2)\right]
\text{d}\overleftarrow{B_s}+\int_t^\infty
2\left( \bar y_j\left( s\right) \right) ^{-}\bar z_j\left( s\right) \text{d}%
\overrightarrow{W_s} \\
&&\ +\int_t^\infty \mathbf{I}_{_{\left\{ \bar y_j(s)\leq 0\right\}
}}\left|
g_j\left( s,y_s^1\right) -g_j\left( s,y_s^2\right) \right| ^2\text{d}s \\
&&\ -\int_t^\infty \mathbf{I}_{\left\{ \bar y_j(s)\leq 0\right\}
}\left| \bar z_s\left( s\right) \right| ^2\text{d}s.
\end{eqnarray*}
\begin{equation}
\tag{3.6}
\end{equation}
From (H4), we have $\xi ^1-\xi ^2\geq 0,$ so
\begin{equation}
\mathbf{E}\left| \left( \xi ^1-\xi ^2\right) ^{-}\right| ^2=0.  \tag{3.7}
\end{equation}
Since $\left( y^1,z^1\right) $ and $\left( y^2,z^2\right) $ are in $\mathcal{%
S}^2\left( [0,\infty ];\mathbf{R}^k\right) \times \mathcal{M}^2\left(
0,\infty ;\mathbf{R}^{k\times d}\right) ,$ it easily follows that
\begin{equation}
\mathbf{E}\int_t^\infty \left( \bar y_j\left( s\right) \right) ^{-}\left[
g_j^1(s,y_s^1)-g_j^1(s,y_s^2)\right] \text{d}\overleftarrow{B_s}=0,
\tag{3.8}
\end{equation}
\begin{equation}
\mathbf{E}\int_t^\infty 2\left( \bar y_j\left( s\right) \right) ^{-}\bar
z_j\left( s\right) \text{d}\overrightarrow{W_s}=0.  \tag{3.9}
\end{equation}
Let
\begin{eqnarray*}
\Delta ^1 &=&-\int_t^\infty 2\left( \bar y_j\left( s\right) \right)
^{-}\left[ f_j^1\left( s,y_s^1,z_s^1\right) -f_j^2\left(
s,y_s^2,z_s^2\right) \right] \text{d}s \\
&=&-\int_t^\infty 2\left( \bar y_j\left( s\right) \right) ^{-}f_j^1\left(
y_1^1,\cdots ,y_j^1,\cdots ,y_k^1,z_1^1,\cdots ,z_j^1,\cdots ,z_k^1\right)
\text{d}s \\
&&+\int_t^\infty 2\left( \bar y_j\left( s\right) \right) ^{-}f_j^1\left(
y_1^1+\bar y_1,\cdots ,y_j^1,\cdots ,y_k^1+\bar y_k^1,z_1^1,\cdots
,z_j^2,\cdots ,z_k^1\right) \text{d}s \\
&&-\int_t^\infty 2\left( \bar y_j\left( s\right) \right) ^{-}f_j^1\left(
y_1^1+\bar y_1,\cdots ,y_j^1,\cdots ,y_k^1+\bar y_k^1,z_1^1,\cdots
,z_j^2,\cdots ,z_k^1\right) \text{d}s \\
&&+\int_t^\infty 2\left( \bar y_j\left( s\right) \right) ^{-}f_j^2\left(
y_1^2,\cdots ,y_j^2,\cdots ,y_k^2,z_1^2,\cdots ,z_j^2,\cdots ,z_k^2\right)
\text{d}s \\
&=&\Delta _1^1+\Delta _2^1,
\end{eqnarray*}
where $\Delta _1^1$ is the first integral, and $\Delta _2^1$ is the second
integral. Since
\[
y_l^1+\bar y_l^{-}\geq y_l^2,\text{ }l\neq j
\]
and (H4), we have
\begin{eqnarray*}
&&f_j^1\left( y_1^1+\bar y_1,\cdots ,y_j^1,\cdots ,y_k^1+\bar
y_k^1,z_1^1,\cdots ,z_j^2,\cdots ,z_k^1\right)  \\
&\geq &f_j^2\left( y_1^2,\cdots ,y_j^2,\cdots ,y_k^2,z_1^2,\cdots
,z_j^2,\cdots ,z_k^2\right) .
\end{eqnarray*}
So $\Delta _2^1\leq 0.$ From (3.6) and the inequality $a^2+b^2\geq 2ab,$ it
follows that
\begin{eqnarray*}
\Delta _1^1 &\leq &2\int_t^\infty \bar y_j\upsilon \left( s\right) \left(
\left| \bar y_1^{-}\right| +\cdots +\left| \bar y_{j-1}^{-}\right| +\left|
\bar y_j\right| +\left| \bar y_{j+1}^{-}\right| +\cdots +\left| \bar
y_k^{-}\right| \right) +\bar y_ju\left( t\right) \left| \bar z_j\right|
\text{d}s \\
&\leq &\int_t^\infty \sum_{i=1}^k\left| \bar y_i^{-}\right| ^2\kappa \left(
s\right) ds+\left( k+d\right) \int_t^\infty \kappa \left( s\right) \left|
\bar y_j\right| ^2\text{d}s \\
&&+\int_t^\infty \sum_{i=1}^d\left| \bar z^{j,i}\right| ^2\mathbf{I}%
_{\left\{ \bar y_j(s)\leq 0\right\} }\text{d}s,
\end{eqnarray*}
where $\kappa \left( t\right) =\upsilon \left( t\right) \vee u^2\left(
t\right) .$ Obviously, $\int_0^\infty \kappa \left( s\right) $d$s<\infty .$
From (H2) we deduce
\begin{eqnarray*}
\Delta ^2 &=&\int_t^\infty \mathbf{I}_{\left\{ \bar y_j(s)\leq 0\right\}
}\left| g_j\left( s,y_s^1\right) -g_j\left( s,y_s^2\right) \right| ^2\text{d}%
s \\
&=&\int_t^\infty \mathbf{I}_{\left\{ \bar y_j(s)\leq 0\right\} }\left|
g_j\left( s,y_1^1,\cdots ,y_j^1,\cdots ,y_k^1\right) -g_j\left(
s,y_1^2,\cdots ,y_j^2,\cdots ,y_k^2\right) \right| ^2\text{d}s \\
&\leq &\int_t^\infty \kappa \left( s\right) \left( \left| \bar
y_1^{-}\right| ^2+\cdots +\left| \bar y_j^{-}\right| ^2+\cdots +\left| \bar
y_k^{-}\right| ^2\right) \text{d}s.
\end{eqnarray*}
So there exists a positive constant $M$, such that
\begin{equation}
\Delta ^1+\Delta ^2\leq M\sum_{i=1}^k\int_t^\infty \kappa \left( s\right)
|\bar y_i^{-}|^2\text{d}s+\int_t^\infty \sum_{i=1}^d|\bar z^{j,i}|^2\mathbf{I%
}_{\left\{ \bar y_j(s)\leq 0\right\} }\text{d}s.  \tag{3.10}
\end{equation}
Taking expectation on both sides of (3.6) and noting (3.7)-(3.10), we get
\[
\mathbf{E}\left| \bar y_j\left( t\right) ^{-}\right| ^2\leq M\mathbf{E}%
\sum_{i=1}^k\int_t^\infty \kappa \left( s\right) \left| \bar y_i^{-}\right|
^2\text{d}s.
\]
Sum from $1$ to $k$ and it follows that
\[
\mathbf{E}\sum_{j=1}^k\left| \bar y_j\left( t\right) ^{-}\right| ^2\leq
Mk\int_t^\infty \kappa \left( s\right) \mathbf{E}\sum_{i=1}^k\left| \bar
y_i^{-}\right| ^2\text{d}s.
\]
By Lemma 5, it follows that
\[
\sum_{j=1}^k\left| \bar y_j(t)^{-}\right| ^2=0,\text{\qquad }\forall t\in
\left[ 0,\infty \right] .
\]
So $\bar y_j\left( t\right) ^{-}=0,$ $j=1,$ $2,\cdots ,$ $k.$ The
proof is completed.\quad $\Box$
 At last, we give an example of
infinite horizon BDSDEs to show our comparison theorem.

\begin{example}
We consider the following BDSDEs,
\begin{equation}
y_t^1=\xi ^1+\int_t^\infty \frac{y_s^1+z_s^1}{1+s^2}\text{d}s+\int_t^\infty
\frac{y_s^1}{1+s^2}\text{d}\overleftarrow{B_s}-\int_t^\infty z_s^1\text{d}%
\overrightarrow{W_s},  \tag{4.1}
\end{equation}
\begin{equation}
y_t^2=\xi ^2+\int_t^\infty \frac{y_s^2+z_s^2}{1+s^2}\text{d}s+\int_t^\infty
\frac{y_s^2}{1+s^2}\text{d}\overleftarrow{B_s}-\int_t^\infty z_s^2\text{d}%
\overrightarrow{W_s},  \tag{4.2}
\end{equation}
where $k=d=l=1,$ $\xi ^i\in L^2\left( \Omega ,\mathcal{F}_\infty ,P;\mathbf{R%
}\right) ,$ $i=1,2$ and $\xi ^1\geq \xi ^2.$ It is easy to check that (4.1)
and (4.2) satisfy (H1) (H2) and (H3). So according to Proposition 3, there
exist unique solutions $\left( y_t^1,z_t^1\right) $ and $\left(
y_t^2,z_t^2\right) $ for (4.1) and (4.2), repectively. Moreover, we can get
the explicit forms as following:
\begin{equation}
y_t^1=\mathbf{E}\left[ \left. \xi ^1\exp \left\{ \int_t^\infty \frac 1{1+s^2}%
\text{d}\overleftarrow{B_s}+\int_t^\infty \frac 1{1+s^2}\text{d}%
\overrightarrow{W_s}\right\} \right| \mathcal{F}_t\right] ,  \tag{4.3}
\end{equation}
\begin{equation}
y_t^2=\mathbf{E}\left[ \left. \xi ^2\exp \left\{ \int_t^\infty \frac 1{1+s^2}%
\text{d}\overleftarrow{B_s}+\int_t^\infty \frac 1{1+s^2}\text{d}%
\overrightarrow{W_s}\right\} \right| \mathcal{F}_t\right] .  \tag{4.4}
\end{equation}

Obviously, by Theorem 6, we also obtain that
\[
y_t^1\geq y_t^2,\qquad \forall t\in \mathbf{R}^{+}.
\]
\end{example}

\textbf{Acknowledgement:} The authors thank the referees for their helpful
suggestions and comments.

\end{document}